\documentclass[conference,9pt]{IEEEtran}
\IEEEoverridecommandlockouts

\def\BibTeX{{\rm B\kern-.05em{\sc i\kern-.025em b}\kern-.08em
    T\kern-.1667em\lower.7ex\hbox{E}\kern-.125emX}}
\usepackage{amsmath,graphicx}

\usepackage[caption=false]{subfig}%
\usepackage{amssymb,amsfonts}
\usepackage{textcomp}
\usepackage{xcolor}

\usepackage{here}
\usepackage{enumerate}
\usepackage{bm}
\usepackage{algorithm}
\usepackage{algpseudocode}
\graphicspath{{Fig/}}
\usepackage{amsthm}
\usepackage{comment}
\usepackage{cases}

\usepackage{mathtools}
\usepackage{bm}
\usepackage{comment}
\usepackage{enumitem}
\usepackage{cite}

\usepackage{stmaryrd}
\usepackage{csvsimple}
\usepackage{multicol}
\usepackage{multirow}
\usepackage{array}

\usepackage{mathabx}
\usepackage{hyperref}
\hypersetup{
    colorlinks=true,
    citecolor=blue,
    linkcolor=red,
    urlcolor=orange,
}

\mathtoolsset{showonlyrefs=true}

\newcommand{\argmin}{\mathop{\mathrm{argmin}}\limits}

\newcommand{\inprod}[2]{{\left\langle #1,#2 \right\rangle}}

\newcommand{\T}{\mathsf{T}}

\newcommand{\prox}[1]{\mathrm{prox}_{#1}}
\newcommand{\moreau}[2]{{}^{#2}#1}

\newcommand{\exR}{\mathbb{R}\cup \{+\infty\}}
\newcommand{\dom}[1]{\mathrm{dom}(#1)}

\newcommand{\norm}[1]{\left\lVert #1 \right\rVert}
\newcommand{\abs}[1]{\left\lvert #1 \right\rvert}

\theoremstyle{plain}%

\theoremstyle{definition}
\newtheorem{theorem}{Theorem}[section]

\newtheorem{problem}[theorem]{Problem}
\newtheorem{definition}[theorem]{Definition}

\newtheorem{assumption}[theorem]{Assumption}

\newtheorem{remark}[theorem]{Remark}

\title{\huge \bf A Proximal Variable Smoothing for Nonsmooth Minimization\\ Involving Weakly Convex Composite with MIMO Application\thanks{This work was supported by JSPS Grants-in-Aid (19H04134, 24K23885).}}

\author{\IEEEauthorblockN{Keita Kume\IEEEauthorrefmark{1},
Isao Yamada\IEEEauthorrefmark{2}}
\IEEEauthorblockA{
  Dept. of Information and Communications Engineering, Institute of Science Tokyo, JAPAN\\
Email: \IEEEauthorrefmark{1}kume@sp.ict.e.titech.ac.jp,
\IEEEauthorrefmark{2}isao@sp.ict.e.titech.ac.jp}}

\begin{document}
\maketitle
\renewcommand{\thefootnote}{\alph{footnote})}
\renewcommand{\appendixname}{Appendix}

\begin{abstract}
  We propose a proximal variable smoothing algorithm for nonsmooth optimization problem with sum of three functions involving weakly convex composite function.
  The proposed algorithm is designed as a time-varying forward-backward splitting algorithm with two steps:
  (i) a time-varying forward step with the gradient of a smoothed surrogate function, designed with the Moreau envelope, of the sum of two functions;
  (ii) the backward step with the proximity operator of the remaining function.
  For the proposed algorithm, we present a convergence analysis in terms of a stationary point by using a newly smoothed surrogate stationarity measure.
  As an application of the target problem, we also present a formulation of multiple-input-multiple-output (MIMO) signal detection with phase-shift keying.
  Numerical experiments demonstrate the efficacy of the proposed formulation and algorithm.
\end{abstract}

\section{Introduction} \label{sec:introduction}
In this paper, we address the following (possibly nonsmooth and nonconvex) optimization problem with sum of three functions:
\begin{problem}\label{problem:origin}
Let
$\mathcal{X}$
and
$\mathcal{Z}$
be Euclidean spaces, i.e., finite-dimensional real Hilbert spaces.
Then,
\begin{equation}
  \mathrm{find} \ \bm{x}^{\star} \in \argmin_{\bm{x}\in \mathcal{X}} (h+g\circ \mathfrak{S} + \phi)(\bm{x}) (\eqqcolon (f + \phi)(\bm{x})),
\end{equation}
where
$\phi$,
$h$,
$\mathfrak{S}$,
$g$,
and
$f\coloneqq h+g\circ\mathfrak{S}$
satisfy the following:
\begin{enumerate}[label=(\roman*)]
  \item
        $\phi:\mathcal{X}\to \exR$
        is a proper lower semicontinuous convex function, and~{\em prox-friendly}, i.e., the proximity operator
        $\prox{\gamma \phi}\ (\gamma \in \mathbb{R}_{++})$
        (see~\eqref{eq:prox})
        has a closed-form expression
        (such a $\phi$ includes, e.g., the indicator function associated with a nonempty (simple) closed convex set);
  \item
        $h:\mathcal{X} \to \mathbb{R}$
        is differentiable and its gradient
        $\nabla h:\mathcal{X}\to\mathcal{X}$
        is Lipschitz continuous over
        $\dom{\phi}\coloneqq\{\bm{x} \in \mathcal{X} \mid \phi(\bm{x}) < +\infty\}$;
  \item
        $\mathfrak{S}: \mathcal{X} \to \mathcal{Z}$
        is a continuously differentiable  mapping;
  \item
        $g:\mathcal{Z} \to \mathbb{R}$
        is (a) Lipschitz continuous (possibly nonsmooth),
        (b)~$\eta$-weakly convex with
        $\eta > 0$,
        i.e.,
        $g+\frac{\eta}{2}\norm{\cdot}^{2}$
        is convex,
        and (c) prox-friendly
        (such a
        $g$
        includes, e.g., $\ell_{1}$-norm
        $\norm{\cdot}_{1}$~\cite{Bauschke-Combettes17}, minimax concave penalty (MCP)~\cite{Zhang10}, and smoothly clipped absolute deviation (SCAD)~\cite{Fan-Li01});
  \item
        $\inf\{(f+\phi)(\bm{x})\mid \bm{x} \in \dom{\phi} \} > -\infty$.
\end{enumerate}
\end{problem}

Problem~\ref{problem:origin} has attracted tremendous attentions because of its potential for flexible formulation to align with many signal processing and machine learning applications (see, e.g., a recent review paper~\cite{Li-So-Ma20}).
Indeed, in a special case where
$\mathfrak{S}$
is linear, such applications include, e.g.,
tensor completion~\cite{Gandy-Recht-Yamada11},
principal component pursuit~\cite{Yin-Parekh-Selesnick19},
and block sparse recovery~\cite{Kuroda-Kitahara22}, to name a few (see, e.g.,~\cite{Combettes-Pesquet21}).
In more general cases where
$\mathfrak{S}$
is nonlinear but continuously differentiable, Problem~\ref{problem:origin} is also known to have a wide range of applications, e.g., robust phase retrieval~\cite{Zhen-Ma-Xue24,Duchi-Ruan18}, robust matrix recovery~\cite{Charisopoulos-Chen-Davis-Diaz-Ding-Drusvyatskiy21,Wang-So-Zoubir23}, and robust blind deconvolution~\cite{Charisopoulos-Davis-Dias-Drusvyatskiy20}, to name a few (see, e.g.,~\cite{Lewis-Wright16,Drusvyatskiy-Paquette19}).

For Problem~\ref{problem:origin} with some special cases, proximal splitting-type algorithms have been used.
For example, in a convex case where
$h$,
$g$,
and
$\phi$
are convex and
$\mathfrak{S}$
is linear, Davis-Yin's three-operator splitting algorithm~\cite{Davis-Yin17,Zhao-Cevher18} has been designed with the gradient of
$h$
and the proximity operators of
$g$
and
$\phi$.
In particular with
$g = 0$,
we can use the forward-backward splitting algorithm (also known as the proximal gradient method), e.g.,~\cite{Bauschke-Combettes17,Condat-Kitahara-Contreras-Hirabayashi23}, which consists of (i) the forward step with the gradient
$\nabla h$
and (ii) the backward step with
$\prox{\gamma \phi}\ (\gamma > 0)$.
For general Problem~\ref{problem:origin}, the so-called prox-linear method~\cite{Lewis-Wright16,Drusvyatskiy-Paquette19} is designed with a solution to a certain proximal subproblem that is a minimization of a linearized function of
$g\circ\mathfrak{S}$.
Although the prox-linear method generates a sequence whose accumulation point is a stationary point in the sense of zero of {\em the general subdifferential} (see~\eqref{eq:general_subdifferential}), the prox-linear method requires some iterative solver for the proximal subproblem at every iteration.

Inspired by~\cite{Bohm-Weight21,Kume-Yamada24,Liu-Xia24,Kume-Yamada24B},
in this paper, we propose for Problem~\ref{problem:origin} a proximal variable smoothing which does not require any iterative solver at each iteration.
Indeed, the proposed algorithm in Alg.~\ref{alg:proposed} is designed
as a time-varying forward-backward splitting-type (or proximal gradient-type) algorithm designed with {\em the Moreau envelope}
$\moreau{g}{\mu}\ (\mu \in (0,\eta^{-1}))$
of
$g$
(see~\eqref{eq:moreau}).
More precisely, the proposed algorithm consists of two steps:
(i) the forward step with the gradient of a time-varying smoothed surrogate function
$f_{n}\coloneqq h+\moreau{g}{\mu_{n}}\circ \mathfrak{S}\ (\mu_{n}(\in (0,\eta^{-1}))\searrow 0)$
of
$f=h+g\circ \mathfrak{S}$
and
(ii) the backward step with
$\prox{\gamma_{n} \phi}$
(see~\eqref{eq:prox})
with a stepsize
$\gamma_{n} > 0$.

We also present an asymptotic convergence analysis, in Theorem~\ref{theorem:convergence_extension}, of the proposed algorithm in terms of a stationary point.
The proposed convergence analysis relies on a smooth approximation of a stationarity measure
$\mathcal{M}_{\gamma}^{f,\phi}\ (\gamma > 0)$~\cite{Liu-Xia24}
in~\eqref{eq:measure} with the smoothed surrogate function
$f_{n}$ (see also Theorem~\ref{theorem:measure}).

The proposed algorithm and analysis serve as an extension of variable smoothing-type algorithms~\cite{Bohm-Weight21,Kume-Yamada24,Liu-Xia24,Kume-Yamada24B} because Problem~\ref{problem:origin} covers the cases assumed in~\cite{Bohm-Weight21,Kume-Yamada24,Liu-Xia24,Kume-Yamada24B} (see Remark~\ref{remark:extension}).
For example, convex constrained minimizations of the composition of a (weakly) convex function and a smooth function have been key models in signal processing (see, e.g.,~\cite{Wang-So-Zoubir23,Charisopoulos-Davis-Dias-Drusvyatskiy20,Lewis-Wright16,Drusvyatskiy-Paquette19}).
However, the existing algorithms~\cite{Bohm-Weight21,Kume-Yamada24,Liu-Xia24,Kume-Yamada24B} cannot be applied to such cases, while the proposed algorithm can be applied.

To verify the effectiveness of the proposed algorithm, we present numerical experiments in a scenario of MIMO signal detection~\cite{Yang-Hanzo15,Chen19,Hayakawa19,Hayakawa-Hayashi20} with its new formulation via Problem~\ref{problem:origin} (see Section~\ref{sec:numerical}).

\vspace{1em}
\noindent
{\bf Notation}:
$\mathbb{N}$,
$\mathbb{Z}$,
$\mathbb{R}$,
$\mathbb{R}_{++}$,
and
$\mathbb{C}$
denote respectively the sets of all positive integers, all integers, all real numbers, all positive real numbers, and all complex numbers
($i$ stands for the imaginary unit, and
$\mathfrak{R}$
and
$\mathfrak{I}$
stand respectively for real and imaginary parts).
The symbol
$\bm{1}$
is the vector of all ones.
For a given matrix
$\bm{A} \in \mathbb{K}^{L\times M}$
and a vector
$\bm{v} \in \mathbb{K}^{N}$
$(\mathbb{K}=\mathbb{R},\mathbb{C})$,
$[\bm{A}]_{l,m} \in \mathbb{K}$
and
$[\bm{v}]_{n} \in \mathbb{K}$
denote respectively the
$(l,m)$
entry of
$\bm{A}$
and
$n$th entry of
$\bm{v}$.
The symbol
$\odot$
is the element-wise product, i.e.,
$[\bm{a}\odot \bm{b}]_{n}\coloneqq [\bm{a}]_{n}[\bm{b}]_{n}\ (\bm{a},\bm{b} \in \mathbb{R}^{N})$,
$\|\cdot\|$
and
$\inprod{\cdot}{\cdot}$
are respectively the Euclidean norm and the standard inner product.
For a linear operator
$A:\mathcal{X}\to \mathcal{Z}$,
$\norm{A}_{\rm op}\coloneqq \sup_{\norm{\bm{x}}\leq 1, \bm{x}\in \mathcal{X}} \norm{A\bm{x}}$
denotes the operator norm of
$A$.

For a continuously differentiable mapping
$\mathcal{F}:\mathcal{X}\to \mathcal{Z}$,
its G\^{a}teaux derivative at
$\bm{x}\in \mathcal{X}$
is the linear operator
$\mathrm{D}\mathcal{F}(\bm{x}):\mathcal{X} \to \mathcal{Z}:\bm{v}\mapsto \lim\limits_{\mathbb{R}\setminus\{0\} \ni t\to 0}\frac{\mathcal{F}(\bm{x}+t\bm{v}) -\mathcal{F}(\bm{x})}{t}$.
For a G\^{a}teaux differentiable function
$J:\mathcal{X} \to \mathbb{R}$,
$\nabla J(\bm{x}) \in \mathcal{X}$
is the gradient of
$J$
at
$\bm{x} \in \mathcal{X}$
if
$\mathrm{D}J(\bm{x})[\bm{v}] = \inprod{\nabla J(\bm{x})}{\bm{v}}\ (\bm{v} \in \mathcal{X})$.
For a function
$J:\mathcal{X} \to \exR$,
$J$
is called (a) proper if
$\dom{J}\coloneqq \{\bm{x} \in \mathcal{X} \mid J(\bm{x}) < +\infty\} \neq \emptyset$;
(b)~lower semicontinuous if
$\{(\bm{x},a) \in \mathcal{X}\times \mathbb{R} \mid J(\bm{x})\leq a\} \subset \mathcal{X}\times \mathbb{R}$
is closed in
$\mathcal{X} \times \mathbb{R}$,
or equivalently,
$\liminf_{n\to\infty}J(\bm{x}_{n}) \geq J(\widebar{\bm{x}})$
for every
$\widebar{\bm{x}}\in \mathcal{X}$
and every
$(\bm{x}_{n})_{n=1}^{\infty} \subset \mathcal{X}$
converging to
$\widebar{\bm{x}}$;
(c) convex if
$J(t\bm{x}_{1}+t\bm{x}_{2}) \leq t J(\bm{x}_{1}) + (1-t)J(\bm{x}_{2})\ (\forall t \in [0,1], \forall \bm{x}_{1},\bm{x}_{2} \in \mathcal{X})$.

\section{Preliminary on Nonsmooth Analysis}
We review some necessary notions in nonsmooth analysis primarily based on notations in~\cite{Rockafellar-Wets98,Bauschke-Combettes17} (see also a review paper~\cite{Li-So-Ma20}).
\begin{definition}[Subdifferential~{\cite[Def. 8.3]{Rockafellar-Wets98}}]\label{definition:subdifferential}
  For a possibly nonconvex function
  $J:\mathcal{X} \to \exR$,
  {\it the general  subdifferential}
  $\partial J:\mathcal{X} \rightrightarrows \mathcal{X}$
  of
  $J$
  at
  $\widebar{\bm{x}}\in \dom{J}$
  is defined by
  \begin{equation}
    \hspace{-0.65em}
    \thickmuskip=0.2\thickmuskip
    \medmuskip=0.2\medmuskip
    \thinmuskip=0.2\thinmuskip
    \arraycolsep=0.2\arraycolsep
    \partial J(\widebar{\bm{x}}) \coloneqq
    \left\{\bm{v} \in \mathcal{X} \mid \substack{\displaystyle\exists (\bm{x}_{n})_{n=1}^{\infty} \subset \mathcal{X},\ \exists \bm{v}_{n} \in \widehat{\partial}J(\bm{x}_{n})\  \mathrm{with}\\\displaystyle  \widebar{\bm{x}}=\lim_{n\to\infty}\bm{x}_{n},\;J(\widebar{\bm{x}}) = \lim_{n\to\infty}J(\bm{x}_{n}), \; \bm{v} = \lim_{n\to\infty}\bm{v}_{n} } \right\}_{\displaystyle,}
    \label{eq:general_subdifferential}
  \end{equation}
  where\footnote{
    The {\em limit inferior}~\cite[Definition 1.5]{Rockafellar-Wets98} is given by
    \begin{equation}
      \textstyle
      \liminf\limits_{\substack{\mathcal{X}\ni \bm{x}\to \widebar{\bm{x}}\\\bm{x}\neq \widebar{\bm{x}}}} \frac{J(\bm{x})-J(\widebar{\bm{x}}) - \inprod{\bm{v}}{\bm{x}-\widebar{\bm{x}}}}{\norm{\bm{x}-\widebar{\bm{x}}}} =
      \sup\limits_{\epsilon>0} \left(\inf\limits_{0<\norm{\bm{x}-\widebar{\bm{x}}}< \epsilon}\frac{J(\bm{x})-J(\widebar{\bm{x}}) - \inprod{\bm{v}}{\bm{x}-\widebar{\bm{x}}}}{\norm{\bm{x}-\widebar{\bm{x}}}}\right).
    \end{equation}
  }
  $\widehat{\partial} J(\widebar{\bm{x}})
    \coloneqq
    \left\{
    \begin{array}{l}
      \substack{\displaystyle                                            \\\displaystyle \\\displaystyle \bm{v} \in \mathcal{X}\mid}
      \liminf\limits_{\substack{\mathcal{X}\ni\bm{x}\to \widebar{\bm{x}} \\ \bm{x}\neq \widebar{\bm{x}}}} \frac{J(\bm{x}) - J(\widebar{\bm{x}}) - \inprod{\bm{v}}{\bm{x}-\widebar{\bm{x}}}}{\norm{\bm{x}-\widebar{\bm{x}}}}  \geq 0
    \end{array}
    \right\}$
  is called {\it the regular  subdifferential} of
  $J$
  at
  $\widebar{\bm{x}} \in \dom{J}$,
  and
  $\partial J(\widebar{\bm{x}})$
  and
  $\widehat{\partial} J(\widebar{\bm{x}})$
  at
  $\widebar{\bm{x}} \notin \dom{J}$
  are understood as
  $\emptyset$.
  If
  $J$
  is convex, then
  $\partial J$
  and
  $\widehat{\partial} J$
  coincide with the convex subdifferential of
  $J$.
\end{definition}

By Fermat's rule~\cite[Theorem 10.1]{Rockafellar-Wets98}, a local minimizer
$\bm{x}^{\star} \in \mathcal{X}$
of
$f+\phi$
in Problem~\ref{problem:origin} satisfies the {\em first-order optimality condition}:
\begin{equation}
  \partial (f+\phi)(\bm{x}^{\star}) (=\partial f(\bm{x}^{\star}) + \partial \phi(\bm{x}^{\star})) \ni \bm{0}, \label{eq:Fermat}
\end{equation}
where the equality can be checked by~\cite[Corollary 10.9]{Rockafellar-Wets98}.
We call
$\bm{x}^{\star}$
satisfying~\eqref{eq:Fermat} {\em stationary point} of
$f+\phi$,
and we focus on finding a stationary point of
$f+\phi$
throughout this paper.

  {\it The proximity operator} and {\it the Moreau envelope} have been used as computational tools for nonsmooth optimization~\cite{Yamada-Yukawa-Yamagishi11,Bauschke-Combettes17,Abe-Yamagishi-Yamada20,Bauschke-Moursi-Wang21}.
For an
$\eta(>0)$-weakly convex function
$g$,
its proximity operator and Moreau envelope of index
$\mu \in (0,\eta^{-1})$
are respectively defined by
\begin{align}
  (\widebar{\bm{z}} \in \mathcal{Z}) \quad & \prox{\mu g}(\widebar{\bm{z}})  \coloneqq \argmin_{\bm{z} \in \mathcal{Z}} \left(g(\bm{z}) + \frac{1}{2\mu}\|\bm{z}-\widebar{\bm{z}}\|^{2}\right); \label{eq:prox} \\
  (\widebar{\bm{z}} \in \mathcal{Z}) \quad & \moreau{g}{\mu}(\widebar{\bm{z}}) \coloneqq \min_{\bm{z} \in \mathcal{Z}} \left(g(\bm{z}) + \frac{1}{2\mu}\|\bm{z}-\widebar{\bm{z}}\|^{2}\right),
  \label{eq:moreau}
\end{align}
where
$\prox{\mu g}$
is single-valued due to the strong convexity of
$g + (2\mu)^{-1}\|\cdot-\widebar{\bm{z}}\|^{2}$.
Many functions
$g$
have closed-form expressions of
$\prox{\mu g}$,
e.g.,
$\ell_{1}$-norm, MCP, and SCAD (see, e.g.,~\cite{prox_repository} and~\cite[Sect. 2.2]{Kume-Yamada24B}).
The Moreau envelope
$\moreau{g}{\mu}$
serves as a smoothed surrogate function of
$g$
because
$\lim_{\mu\to 0}\moreau{g}{\mu}(\bm{z}) = g(\bm{z})\ (\bm{z} \in \mathcal{Z})$,
and
$\moreau{g}{\mu}$
is continuously differentiable with
$\nabla \moreau{g}{\mu}(\bm{z}) = \mu^{-1}(\bm{z}-\prox{\mu g}(\bm{z}))$~\cite{Bauschke-Moursi-Wang21}.
Moreover,
$\nabla \moreau{g}{\mu}$
is Lipschitz continuous~\cite[Corollary 3.4]{Hoheisel-Laborde-Oberman20}.

\section{Proximal variable smoothing}
We employ the following stationarity measure~\cite{Liu-Xia24} at a given
$\widebar{\bm{x}} \in \mathcal{X}$
with
$\gamma \in \mathbb{R}_{++}$
in order to measure an achievement level of the first-order optimality condition in~\eqref{eq:Fermat}:
\begin{equation}
  \mathcal{M}_{\gamma}^{f,\phi}(\widebar{\bm{x}})
  \coloneqq \inf\left\{\gamma^{-1}\norm{\widebar{\bm{x}} - \prox{\gamma\phi}(\widebar{\bm{x}}-\gamma \bm{v})} \mid \bm{v} \in \partial f(\widebar{\bm{x}})\right\}. \label{eq:measure}
\end{equation}
The measure
$\mathcal{M}_{\gamma}^{f,\phi}$
serves as a generalization of commonly-used stationarity measures.
Indeed, with a special case
$\phi \equiv 0$,
$\mathcal{M}_{\gamma}^{f,\phi}(\widebar{\bm{x}}) = d(\bm{0},\partial f(\widebar{\bm{x}}))\ (\widebar{\bm{x}}\in \mathcal{X})$
has been used as a standard stationarity measure in nonsmooth optimization (see, e.g.,~\cite{Rockafellar-Wets98,Bohm-Weight21,Kume-Yamada24,Kume-Yamada24B}), where
$d(\cdot,\cdot)$
stands for the distance between a given point and a given subset in
$\mathcal{X}$.
Even for
$\phi \neq 0$,
$\mathcal{M}_{\gamma}^{f,\phi}$
works as a stationarity measure at
$\widebar{\bm{x}} \in \mathcal{X}$
for Problem~\ref{problem:origin}
because~\cite[p.243]{Liu-Xia24}
\begin{equation}
  (\forall \gamma \in \mathbb{R}_{++}) \quad \partial f(\widebar{\bm{x}})+\partial\phi(\widebar{\bm{x}})  \ni \bm{0}
  \Leftrightarrow \mathcal{M}_{\gamma}^{f,\phi}(\widebar{\bm{x}})=0. \label{eq:Fermat_measure}
\end{equation}
Thus, finding a stationary point
$\bm{x}^{\star}\in\mathcal{X}$
of Problem~\ref{problem:origin} is equivalent to finding
$\bm{x}^{\star}$
such that
$\mathcal{M}_{\gamma}^{f,\phi}(\bm{x}^{\star}) = 0$.

From a viewpoint of approximating a stationary point
$\bm{x}^{\star} \in \mathcal{X}$
with an iterative update of estimates
$(\bm{x}_{n})_{n=1}^{\infty} \subset \mathcal{X}$
of
$\bm{x}^{\star}$,
a stationarity measure is desired to have a lower semicontinuous property\footnote{
  To explain this, let
  $\mathcal{M}:\mathcal{X} \to \mathbb{R}$
  be a stationarity measure without lower semicontinuity.
  In this case,
  even if any cluster point of
  $(\bm{x}_{n})_{n=1}^{\infty} \subset \mathcal{X}$
  is not a stationary point, then
  $\liminf_{n\to\infty}\mathcal{M}(\bm{x}_{n})=0$
  may hold~\cite[p. 837]{Levin-Kileel-Boumal23}.
}.
Fortunately,
$\mathcal{M}_{\gamma}^{f,\phi}$
in~\eqref{eq:measure} of our interest is lower semicontinuous, and
$\mathcal{M}_{\gamma}^{f,\phi}$
can be approximated with a smoothed surrogate function
$h+\moreau{g}{\mu}\circ\mathfrak{S}$
of
$f=h+g\circ\mathfrak{S}$.

\begin{theorem} \label{theorem:measure}
  Consider Problem~\ref{problem:origin}.
  For arbitrarily given
  $\widebar{\bm{x}} \in \mathcal{X}$
  and
  $(\bm{x}_{n})_{n=1}^{\infty} \subset\mathcal{X}$
  such that
  $\lim_{n\to\infty}\bm{x}_{n}=\widebar{\bm{x}}$,
  the following hold:
  \begin{enumerate}[label=(\alph*)]
    \item (Lower semicontinuity of $\mathcal{M}_{\gamma}^{f,\phi}$) \label{enum:measure_lsc}
          \begin{equation}
            \hspace{-1.3em}
            (\gamma \in \mathbb{R}_{++}) \ \liminf_{n\to \infty} \mathcal{M}_{\gamma}^{f,\phi}(\bm{x}_{n})\coloneqq \sup_{n\in \mathbb{N}} \inf_{k\geq n}\mathcal{M}_{\gamma}^{f,\phi}(\bm{x}_{k}) \geq \mathcal{M}_{\gamma}^{f,\phi}(\widebar{\bm{x}}).
          \end{equation}
    \item (Smooth approximation of $\mathcal{M}_{\gamma}^{f,\phi}$) \label{enum:measure_smoothed}\label{enum:stationarity}
          Let
          $(\mu_{n})_{n=1}^{\infty} \subset (0,\eta^{-1})$
          satisfy
          $\mu_{n}\searrow 0\ (n\to\infty)$,
          and
          $f_{n}\coloneqq h+\moreau{g}{\mu_{n}}\circ \mathfrak{S}$.
          Then,
          \begin{equation}
            (\gamma \in \mathbb{R}_{++}) \quad \liminf_{n\to\infty} \mathcal{M}_{\gamma}^{f_{n}, \phi}(\bm{x}_{n}) \geq  \mathcal{M}_{\gamma}^{f,\phi}(\widebar{\bm{x}}).
          \end{equation}
          Moreover,
          if
          $\liminf\limits_{n\to\infty} \mathcal{M}_{\gamma}^{f_{n},\phi}(\bm{x}_{n}) = 0$
          with some
          $\gamma \in \mathbb{R}_{++}$,
          then
          $\widebar{\bm{x}}$
          is a stationary point of
          $f+\phi$.
  \end{enumerate}
\end{theorem}

Theorem~\ref{theorem:measure}~\ref{enum:stationarity} implies that our goal for finding a stationary point
$\bm{x}^{\star} \in \mathcal{X}$
of Problem~\ref{problem:origin} is reduced to the following problem:
\begin{equation}
  \thickmuskip=0.3\thickmuskip
  \medmuskip=0.3\medmuskip
  \thinmuskip=0.3\thinmuskip
  \arraycolsep=0.3\arraycolsep
  \mathrm{find\ a\ convergent}\
  (\bm{x}_{n})_{n=1}^{\infty}\subset \mathcal{X}
  \ \mathrm{such\ that}\
  \liminf\limits_{n\to\infty}\mathcal{M}_{\gamma}^{f_{n},\phi}(\bm{x}_{n}) = 0
\end{equation}
with some
$\gamma \in \mathbb{R}_{++}$
and
$f_{n}\coloneqq h+\moreau{g}{\mu_{n}}\circ \mathfrak{S}\ (\mu_{n}\searrow 0)$.
In order to find such a sequence
$(\bm{x}_{n})_{n=1}^{\infty}$,
we propose a proximal variable smoothing illustrated in Algorithm~\ref{alg:proposed}.

\begin{algorithm}[t]
  \caption{Proximal variable smoothing for Problem~\ref{problem:origin}}
  \label{alg:proposed}
  \algrenewcommand\algorithmicindent{0.9em}%
  {
    \begin{algorithmic}[0]
      \Require
      $\bm{x}_{1}\in \dom{\phi},c\in(0,1),(\mu_{n})_{n=1}^{\infty} \subset (0,\frac{1}{2\eta}]$ satisfying~\eqref{eq:nonsummable}
      \For{$n=1,2,\ldots$}
      \State
      Set
      $f_{n}\coloneqq h+\moreau{g}{\mu_{n}}\circ \mathfrak{S}$
      \State
      Find
      $\gamma_{n} > 0$
      satisfying Assumption~\ref{assumption:Lipschitz}~\ref{enum:stepsize}
      {\small (see Remark~\ref{remark:example}~\ref{enum:ex:stepsize})}
      \State
      $\bm{x}_{n+1} \leftarrow \prox{\gamma_{n}\phi}(\bm{x}_{n}-\gamma_{n}\nabla f_{n}(\bm{x}_{n}))$
      \EndFor
      \Ensure
      $\bm{x}_{n} \in \dom{\phi}$
    \end{algorithmic}
  }
\end{algorithm}

The proposed Algorithm~\ref{alg:proposed} is designed as a time-varying forward-backward splitting-type algorithm.
In the forward step at
$n$th iteration,
we perform a gradient descent
at the latest estimate
$\bm{x}_{n} \in \mathcal{X}$
as
$\bm{x}_{n+\frac{1}{2}}\coloneqq \bm{x}_{n} - \gamma_{n} \nabla f_{n}(\bm{x}_{n}) \in \mathcal{X}$
with a smoothed surrogate function
$f_{n}$
and
a stepsize
$\gamma_{n}$.
In the backward step at
$n$th iteration, we assign
$\prox{\gamma_{n}\phi}(\bm{x}_{n+\frac{1}{2}})$
to the next estimate
$\bm{x}_{n+1} \in \mathcal{X}$.

The index
$\mu_{n}\in (0,2^{-1}\eta^{-1}]$
of
$\moreau{g}{\mu_{n}}$
is designed to satisfy
\begin{equation}
  \begin{cases}
     & {\rm(i)}\ \lim_{n\to\infty} \mu_{n} = 0,
    {\rm(ii)}\ \sum_{n=1}^{\infty} \mu_{n} = +\infty,                                                                    \\
     & \displaystyle {\rm(iii)}\ (\exists M \geq 1, \forall n \in \mathbb{N}) \quad M^{-1}\leq \mu_{n+1}/\mu_{n} \leq 1.
  \end{cases}
  \label{eq:nonsummable}
\end{equation}
For example,
$(\mu_{n})_{n=1}^{\infty}\coloneqq((2\eta)^{-1}n^{-1/\alpha})_{n=1}^{\infty}$
with
$\alpha \geq 1$
enjoys~\eqref{eq:nonsummable},
and
$\alpha = 3$
has been used as a standard choice of
$(\mu_{n})_{n=1}^{\infty}$
to achieve a reasonable convergence rate in variable smoothing~\cite{Bohm-Weight21,Liu-Xia24}.

Every stepsize
$\gamma_{n} \in \mathbb{R}_{++}$
in Algorithm~\ref{alg:proposed} is chosen to satisfy an Armijo-type condition (see, e.g.,~\cite{Beck17}) with
$\gamma\coloneqq \gamma_{n}$
and
$c \in (0,1)$:
\begin{equation} \label{eq:Armijo}
  \begin{array}{ll}
     & (f_{n}+\phi)\left(\prox{\gamma\phi}(\bm{x}_{n} - \gamma \nabla f_{n}(\bm{x}_{n}))\right)                             \\
     & \hspace{5em} \leq (f_{n}+\phi)(\bm{x}_{n}) - c\gamma \left(\mathcal{M}_{\gamma}^{f_{n},\phi}(\bm{x}_{n})\right)^{2}.
  \end{array}
\end{equation}
The existence of such a stepsize
$\gamma_{n}$
is guaranteed under the Lipschitz continuity of
$\nabla f_{n}$
(see Assumption~\ref{assumption:Lipschitz} and Remark~\ref{remark:example}).
\begin{assumption}\label{assumption:Lipschitz}
  Consider Problem~\ref{problem:origin} and Algorithm~\ref{alg:proposed}.
  For
  $f_{n} =  h+\moreau{g}{\mu_{n}}\circ\mathfrak{S}$,
  we assume:
  \begin{enumerate}[label=(\alph*)]
    \item (Gradient Lipschitz continuity condition) \label{enum:gradient_Lipschitz}
          There exist
          $\varpi_{1},\varpi_{2} \in \mathbb{R}_{++}$
          such that
          $\nabla f_{n}\ (n\in \mathbb{N})$
          is Lipschitz continuous with a Lipschitz constant
          $L_{\nabla f_{n}}\coloneqq \varpi_{1}+\varpi_{2}\mu_{n}^{-1}.$
    \item (Lower bound condition for stepsizes) \label{enum:stepsize}
          For some
          $c \in (0,1)$
          and
          $\beta \in \mathbb{R}_{++}$,
          $\gamma_{n}\ (n\in \mathbb{N})$
          satisfies (i) the Armijo-type condition~\eqref{eq:Armijo}
          with
          $\gamma\coloneqq \gamma_{n}$;
          (ii)
          $\gamma_{n} \geq \beta L_{\nabla f_{n}}^{-1}$;
          and (iii)
          $\widebar{\gamma}\coloneqq \sup\limits_{n\in \mathbb{N}} \gamma_{n} < +\infty$.
  \end{enumerate}
\end{assumption}

\begin{remark}[Examples achieving Assumption~\ref{assumption:Lipschitz}]\label{remark:example}
  \
  \begin{enumerate}[label=(\alph*)]
    \item \label{enum:Lipschitz_S}
          Assumption~\ref{assumption:Lipschitz}~\ref{enum:gradient_Lipschitz} is satisfied if
          (i)
          $\norm{\mathrm{D}\mathfrak{S}(\cdot)}_{\rm op}$
          is bounded over
          $\dom{\phi}$
          and
          (ii)
          $\mathrm{D}\mathfrak{S}(\cdot)$
          is Lipschitz continuous over
          $\dom{\phi}$ (see, e.g.,~\cite[Prop. 4.5]{Kume-Yamada24B}).
          Clearly, Assumption~\ref{assumption:Lipschitz}~\ref{enum:gradient_Lipschitz} is automatically satisfied if
          $\mathfrak{S}$
          is linear.
    \item \label{enum:ex:stepsize}
          We have mainly two choices (see, e.g.,~\cite{Beck17}) of
          $\gamma_{n}$
          enjoying Assumption~\ref{assumption:Lipschitz}~\ref{enum:stepsize}.
          The first one is
          $\gamma_{n}\coloneqq 2(1-c)L_{\nabla f_{n}}^{-1}$.
          The second one is given by the so-called {\em backtracking algorithm},
          i.e.,
          $\gamma_{n}\coloneqq   \max\{\gamma_{\rm initial}\rho^{k} \mid k\in \mathbb{N}\cup \{0\}, \gamma \coloneqq \gamma_{\rm initial}\rho^{k}\ \mathrm{satisfies~\eqref{eq:Armijo}}\}$
          with
          $\gamma_{\rm initial} \in \mathbb{R}_{++}$
          and
          $\rho \in (0,1)$,
          where such a
          $\gamma_{n}$
          can be obtained in finite arithmetic operations.
          For the second choice, any knowledge on
          $L_{\nabla f_{n}}$
          is not required.
  \end{enumerate}
\end{remark}

We present below a convergence analysis of Algorithm~\ref{alg:proposed}.
\begin{theorem}[Convergence analysis of Alg.~\ref{alg:proposed}]\label{theorem:convergence_extension}
  Consider Problem~\ref{problem:origin}.
  Choose arbitrarily
  $\bm{x}_{1} \in \dom{\phi}$,
  $c \in (0,1)$,
  and
  $(\mu_{n})_{n=1}^{\infty} \subset (0,(2\eta)^{-1}]$
  satisfying~\eqref{eq:nonsummable}.
  Suppose that
  $(\bm{x}_{n})_{n=1}^{\infty} \subset \mathcal{X}$
  is generated by Alg.~\ref{alg:proposed} under Assumption~\ref{assumption:Lipschitz}.
  Then, the following hold:
  \begin{enumerate}[label=(\alph*)]
    \item \label{enum:liminf}
          With
          $\widebar{\gamma} = \sup_{n\in\mathbb{N}}\gamma_{n} < + \infty$
          in Assumption~\ref{assumption:Lipschitz}~\ref{enum:stepsize},
          \begin{equation}
            \liminf_{n\to\infty} \mathcal{M}_{\widebar{\gamma}}^{f_{n},\phi}(\bm{x}_{n}) = 0.
          \end{equation}
    \item
          We can choose a subsequence
          $(\bm{x}_{m(l)})_{l=1}^{\infty} \subset \mathcal{X}$
          of
          $(\bm{x}_{n})_{n=1}^{\infty}$
          such that
          $\lim_{l\to\infty} \mathcal{M}_{\widebar{\gamma}}^{f_{m(l)},\phi}(\bm{x}_{m(l)}) = 0$,
          where
          $m:\mathbb{N}\to\mathbb{N}$
          is monotonically increasing.
          Moreover, every cluster point
          $\bm{x}^{\star} \in \mathcal{X}$
          of
          $(\bm{x}_{m(l)})_{l=1}^{\infty}$
          is a stationary point of Problem~\ref{problem:origin}
          (This statement can be checked by combining~\ref{enum:liminf} with Theorem~\ref{theorem:measure}~\ref{enum:stationarity}).
  \end{enumerate}
\end{theorem}

\begin{remark}[Relation to previous works~\cite{Bohm-Weight21,Kume-Yamada24,Liu-Xia24,Kume-Yamada24B}] \label{remark:extension}
  \
  \begin{enumerate}[label=(\alph*)]
    \item (Variable smoothing in~\cite{Bohm-Weight21,Kume-Yamada24,Kume-Yamada24B})
          For Problem~\ref{problem:origin} with a special case
          $\phi \coloneqq 0$,
          Algorithm~\ref{alg:proposed} reproduces a variable smoothing in~\cite{Bohm-Weight21,Kume-Yamada24,Kume-Yamada24B}.
          Indeed, a stationarity measure used in~\cite{Bohm-Weight21,Kume-Yamada24,Kume-Yamada24B} can be expressed as
          $\mathcal{M}_{\gamma}^{f,\phi}$
          with
          $\phi=0$.
    \item (Proximal variable smoothing in~\cite{Liu-Xia24})
          For Problem~\ref{problem:origin} with a special case where
          $\mathfrak{S}$
          is surjective and linear,
          Algorithm~\ref{alg:proposed} reproduces a proximal variable smoothing in~\cite[Alg. 1]{Liu-Xia24}.
          A convergence analysis in~\cite[Theorem 1]{Liu-Xia24} of~\cite[Alg. 1]{Liu-Xia24}
          provides only a bound of the number of iterations to achieve
          $\mathcal{M}_{\gamma}^{f,\phi}(\bm{x}_{n}) < \epsilon$
          for a given
          $\epsilon \in \mathbb{R}_{++}$.
          In contrast, thanks to Theorem~\ref{theorem:measure},
          Theorem~\ref{theorem:convergence_extension} ensures an asymptotic convergence of Alg.~\ref{alg:proposed} in more general settings than the surjective and linearity of
          $\mathfrak{S}$
          assumed in~\cite{Liu-Xia24}
          (Note: Assumption~\ref{assumption:Lipschitz}~\ref{enum:gradient_Lipschitz} is automatically satisfied if
          $\mathfrak{S}$
          is linear [see Remark~\ref{remark:example}~\ref{enum:Lipschitz_S}]).
  \end{enumerate}
\end{remark}

\section{Application to MU-MIMO signal detection}\label{sec:numerical}
\begin{figure*}[t]
  \centering
  \begin{minipage}[b]{0.55\columnwidth}
    \centering
    \includegraphics[clip,width=\columnwidth]{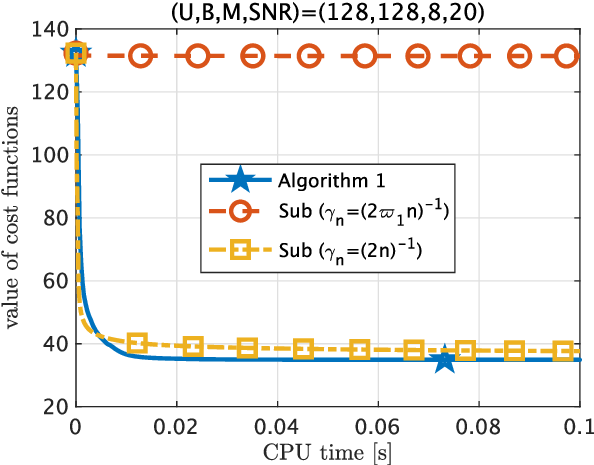}
    \caption{Convergence history}
    \label{fig:convergence}
    \vspace{-1em}
  \end{minipage}
  \hspace{0.02\columnwidth} %
  \begin{minipage}[b]{0.55\columnwidth}
    \centering
    \includegraphics[clip, width=\columnwidth]{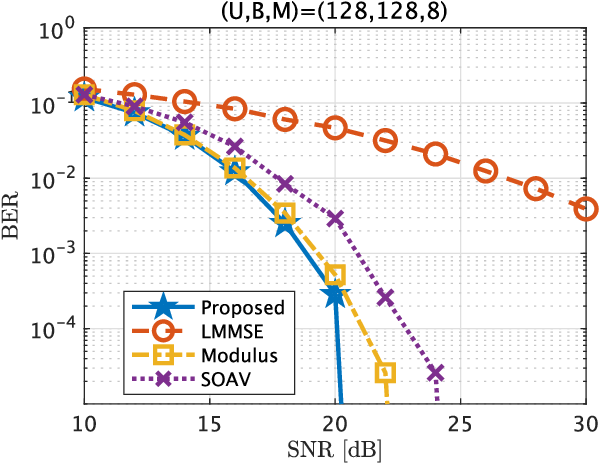}
    \caption{BER vs SNR ($B=U$)}
    \label{fig:BER1}
    \vspace{-1em}
  \end{minipage}
  \hspace{0.02\columnwidth} %
  \begin{minipage}[b]{0.55\columnwidth}
    \centering
    \includegraphics[clip,width=\columnwidth]{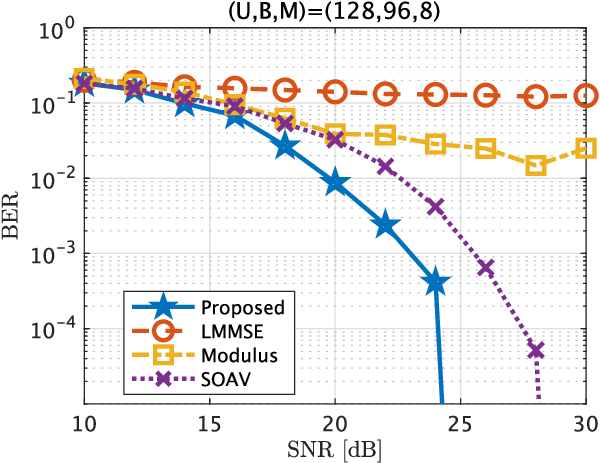}
    \caption{BER vs SNR ($B=\frac{3}{4}U$)}
    \label{fig:BER2}
    \vspace{-1em}
  \end{minipage}
\end{figure*}

\subsection{Formulation of MU-MIMO Signal Detection via Problem~\ref{problem:origin}} \label{sec:numerical_formulation}
To evaluate the performance of Algorithm~\ref{alg:proposed}, we carried out numerical experiments in a scenario of
massive multiuser multiple-input-multiple-output (MU-MIMO) signal detection using
$M(\in \mathbb{N})$-ary
phase-shift keying
(see, e.g.,~\cite{Yang-Hanzo15,Hayakawa-Hayashi20,Chen19,Hayakawa19}).
This task is reduced to
\begin{equation}
  \mathrm{find}\ \mathsf{s}^{\star} \in \mathsf{D} \subset \mathbb{C}^{U}
  \ \mathrm{from}\
  \mathsf{y} = \mathsf{H}\mathsf{s}^{\star} + \mathsf{e} \in \mathbb{C}^{B}, \label{eq:MIMO_model}
\end{equation}
where
$\mathsf{s}^{\star} \in \mathsf{D}\coloneqq \{ \exp(i2\pi m/M) \in \mathbb{C} \mid m = 0, 1,\ldots,M-1\}^{U} \subset \mathbb{C}^{U}$
is the transmit vector,
$\mathsf{y} \in \mathbb{C}^{B}$
is the received vector with a known measurement matrix
$\mathsf{H} \in \mathbb{C}^{B\times U}$,
$\mathsf{e} \in \mathbb{C}^{B}$
is noise,
$\mathsf{D}$
is a discrete set called a {\em constellation set},
and
$U,B\in \mathbb{N}$
are respectively the numbers of transmitting and receiving antennas.
The task~\eqref{eq:MIMO_model} with the case
$U > B$
has appeared in, e.g., typical IoT environments~\cite{Chen19}.

Usually, MU-MIMO signal detection problem in~\eqref{eq:MIMO_model} has been tackled via real-valued optimization problems~\cite{Yang-Hanzo15,Chen19,Hayakawa-Hayashi20} as
\begin{equation}
  \mathop{\mathrm{minimize}}\limits_{\bm{s} \in C \subset \mathbb{R}^{2U}}\ \frac{1}{2}\norm{\bm{y}-\bm{H}\bm{s}}^{2} + \lambda \psi(\bm{s})
  \label{eq:MIMO_problem}
\end{equation}
with the real-valued expressions of
$\mathsf{y}$
and
$\mathsf{H}$:
\begin{equation}
  \hspace{-1em}
  \bm{y}\coloneqq \widehat{\mathsf{y}} \coloneqq \begin{bmatrix} \mathfrak{R}(\mathsf{y})\\ \mathfrak{I}(\mathsf{y}) \end{bmatrix} \in \mathbb{R}^{2B},\;
  \bm{H} \coloneqq \begin{bmatrix} \mathfrak{R}(\mathsf{H}) & - \mathfrak{I}(\mathsf{H}) \\ \mathfrak{I}(\mathsf{H}) & \mathfrak{R}(\mathsf{H}) \end{bmatrix} \in \mathbb{R}^{2B\times 2U},
  \label{eq:real_valued}
\end{equation}
where a constraint set
$(\emptyset \neq )C \subset \mathbb{R}^{2U}$
is designed  to contain the real-valued expression
$\mathfrak{D}\coloneqq \{\widehat{\mathsf{s}} \mid \mathsf{s} \in \mathsf{D}\} \subset \mathbb{R}^{2U}$
of
$\mathsf{D}$,
$\lambda \in \mathbb{R}_{++}$
is a weight
and
$\psi:\mathbb{R}^{2U}\to \mathbb{R}$
is a regularizer.
For example, a classical linear minimum mean-square-error detection (denoted by ``LMMSE''), e.g.,~\cite{Yang-Hanzo15}, considers the problem~\eqref{eq:MIMO_problem} with
$\psi\coloneqq \norm{\cdot}^{2}$
and
$C \coloneqq \mathbb{R}^{2U}$,
where
$\lambda \coloneqq \sigma^{2}/2$
is given by the variance
$\sigma^{2}$
of
$\mathsf{e}$
in~\eqref{eq:MIMO_model}.
In this case, the closed-form solution is given by
$(\bm{H}^{\T}\bm{H}+\sigma^{2}\bm{I})^{-1}\bm{H}^{\T}\bm{y}$.
The reference~\cite{Chen19} proposed a modulus-constrained least squares model~\eqref{eq:MIMO_problem} (denoted by ``Modulus'') where
$\psi \coloneqq 0$,
and
$C \coloneqq \{\widehat{\mathsf{z}} \mid \mathsf{z} \in \mathbb{T}^{U}\} \subset \mathbb{R}^{2U}$
with
$\mathbb{T}\coloneqq \{z \in \mathbb{C} \mid \abs{z} = 1\}$
(Note: the modulus constraint $C$
is designed as a superset of
$\mathfrak{D}$).
Moreover,~\cite{Hayakawa-Hayashi20} proposed the so-called SOAV model~\eqref{eq:MIMO_problem} (denoted by ``SOAV'') where
$C$
is the convex hull of
$\mathfrak{D}$
and
$\psi$
is the {\em sum-of-absolute-value (SOAV)} function~\cite{Nagahara15}.
The SOAV function
$\psi_{\rm SOAV}\coloneqq\frac{1}{M}\sum_{m=1}^{M} \norm{\cdot- \widehat{\mathsf{s}}_{m}}_{1}$
is designed to penalize a given
$\bm{s} \in \mathbb{R}^{2U}$
that deviates from discrete-valued points in
$\mathfrak{D}$,
where
$\mathsf{s}_{m} \coloneqq \exp(i2\pi m/M)\bm{1} \in \mathsf{D}$
and
$\widehat{\mathsf{s}}_{m} \in \mathfrak{D}\ (m=0,1,\ldots,M-1)$
(see~\eqref{eq:real_valued} for $\widehat{(\cdot)}$).
However, the penalization by
$\psi_{\rm SOAV}$
is not contrastive enough to distinguish between points in
$\mathfrak{D}$
and the other points because any point in
$\mathfrak{D}$
is never unique minimizer of
$\psi_{\rm SOAV}$
due to its convexity\footnote{
As more contrastive regularizers than the SOAV function
$\psi_{\rm SOAV}$,
nonconvex regularizers
$\psi\coloneqq \frac{1}{M}\sum_{m=1}^{M} \norm{\cdot-\widehat{\mathsf{s}}_{m}}_{p}$
have been proposed~\cite{Hayakawa19} with
$\ell_{p}$-pseudonorm
$\norm{\cdot}_{p}\ (p\in[0,1))$.
However, for the problem~\eqref{eq:MIMO_problem} with this nonconvex
$\psi$,
any solver with guaranteed convergence has not been found yet as remarked in~\cite{Hayakawa19}
due to the severe nonconvexity.
}.

In this paper, for~\eqref{eq:MIMO_model}, we propose a new model with a more contrastive regularizer
$\psi_{\lambda_{r},\lambda_{\theta}}$
than the SOAV function
$\psi_{\rm SOAV}$
as:
\begin{equation}
  \hspace{-1.5em}
  \thickmuskip=0.0\thickmuskip
  \medmuskip=0.0\medmuskip
  \thinmuskip=0.0\thinmuskip
  \arraycolsep=0.0\arraycolsep
  \mathop{\mathrm{minimize}}\limits_{\substack{\bm{r} \in [\underline{r},1]^{U}_{,} \\ \bm{\theta}\in \mathbb{R}^{U}}}
  \frac{1}{2}\norm{\bm{y}-\bm{H}F(\bm{r},\bm{\theta})}^{2}
  + \underbrace{\lambda_{r}\sum_{u=1}^{U} [\bm{r}]_{u}^{-1}
  + \lambda_{\theta}\norm{\mathrm{\bf sin}\left(\frac{M\bm{\theta}}{2}\right)}_{1}}_{\eqqcolon \psi_{\lambda_{r},\lambda_{\theta}}(\bm{r},\bm{\theta})}\hspace{-0.03em},
  \label{eq:MIMO_proposed}
\end{equation}
where we use a polar coordinate-type expression of
$\bm{s}$
in~\eqref{eq:MIMO_problem} as:
\begin{equation}
  F:\mathbb{R}^{U}\times \mathbb{R}^{U} \to \mathbb{R}^{2U}:(\bm{r},\bm{\theta}) \mapsto \bm{s}\coloneqq
  \begin{bmatrix}
    \bm{r}\odot \mathrm{\bf sin}(\bm{\theta}) \\ \bm{r} \odot \mathrm{\bf cos}(\bm{\theta})
  \end{bmatrix},
\end{equation}
i.e.,
$\bm{r}$
and
$\bm{\theta}$
denote respectively the radial and angle coordinates of
$\bm{s}$,
{
\thickmuskip=0.6\thickmuskip
\medmuskip=0.6\medmuskip
\thinmuskip=0.6\thinmuskip
\arraycolsep=0.6\arraycolsep
$\mathrm{\bf sin}:\mathbb{R}^{U}\times\mathbb{R}^{U}\to \mathbb{R}^{U}:\bm{\theta}\mapsto \begin{bmatrix} \sin([\bm{\theta}]_{1}) & \sin([\bm{\theta}]_{2}) & \cdots & \sin([\bm{\theta}]_{U}) \end{bmatrix}^{\T}$
}
(so as $\mathrm{\bf cos}$),
$\underline{r} \in (0,1]$
is the smallest value for each
$[\bm{r}]_{u}$,
and
$\lambda_{r},\lambda_{\theta} \in \mathbb{R}_{++}$
are weights.

The proposed
$\psi_{\lambda_{r},\lambda_{\theta}}$
in~\eqref{eq:MIMO_proposed}
penalizes
(i) small radius
$[\bm{r}]_{u}$
and
(ii)
angle
$[\bm{\theta}]_{u}$
that deviates from desired angles
$\{ 2\pi m/M \mid m \in \mathbb{Z}\}$
(see
$\mathsf{D}$
just after~\eqref{eq:MIMO_model}).
Indeed,
for any
$\lambda_{r},\lambda_{\theta} \in \mathbb{R}_{++}$,
we can check:
\begin{align}
  \hspace{-1em}
  F(\bm{r}^{\star},\bm{\theta}^{\star}) \in \mathfrak{D} \Leftrightarrow (\bm{r}^{\star},\bm{\theta}^{\star}) \in & \argmin_{(\bm{r},\bm{\theta}) \in [\underline{r},1]^{U}\times \mathbb{R}^{U}} \psi_{\lambda_{r},\lambda_{\theta}}(\bm{r},\bm{\theta}) \\
                                                                                                                  & \ (=\{(\bm{1}, 2\pi \bm{m}/M) \mid \bm{m} \in \mathbb{Z}^{U}\}), \label{eq:desired}
\end{align}
implying thus
$\psi_{\lambda_{r},\lambda_{\theta}}$
is a desired regularizer to penalize a point
$\bm{s} = F(\bm{r},\bm{\theta})$
that deviates from points in
$\mathfrak{D}$.
To examine a potential of the basic idea\footnote{
  The desired property~\eqref{eq:desired} holds even if
  $\norm{\cdot}_{1}$
  in
  $\psi_{\lambda_{r},\lambda_{\theta}}$
  is replaced by weakly convex functions, e.g.,
  MCP~\cite{Zhang10}, and SCAD~\cite{Fan-Li01}.
} in
$\psi_{\lambda_{r},\lambda_{\theta}}$,
we will present the numerical performance of the model~\eqref{eq:MIMO_proposed}, which
can be formulated as Problem~\ref{problem:origin}\footnote{
$g$
and
$\phi$
are prox-friendly (see, e.g.,~\cite[Exm. 24.20]{Bauschke-Combettes17} and~\cite[Lem. 6.26]{Beck17}).
} with
$\mathcal{X}\coloneqq \mathbb{R}^{U}\times \mathbb{R}^{U}$,
$\mathcal{Z} \coloneqq \mathbb{R}^{U}$,
and
{
  \thickmuskip=0.3\thickmuskip
  \medmuskip=0.3\medmuskip
  \thinmuskip=0.3\thinmuskip
  \arraycolsep=0.3\arraycolsep
\begin{align}
   & h:\mathcal{X}\to\mathbb{R}: (\bm{r},\bm{\theta}) \mapsto \frac{1}{2}\norm{\bm{y}-\bm{H}F(\bm{r},\bm{\theta})}^{2} + \lambda_{r}\sum_{u=1}^{U}[\bm{r}]_{u}^{-1},                             \\
   & \mathfrak{S}:\mathcal{X}\to\mathcal{Z}:(\bm{r},\bm{\theta}) \mapsto \mathrm{\bf sin}\left(\frac{M\bm{\theta}}{2}\right),
  \quad
  g:\mathcal{Z}\to\mathbb{R}:\bm{z} \mapsto \lambda_{\theta}\norm{\bm{z}}_{1},                                                                                                                         \\
   & \phi:\mathcal{X}\to\exR:(\bm{r},\bm{\theta}) \mapsto \begin{cases} 0,       & \mathrm{if}\ (\bm{r},\bm{\theta})\in  [\underline{r},1]^{U}\times \mathbb{R}^{U};    \\
              +\infty, & \mathrm{if}\ (\bm{r},\bm{\theta}) \notin [\underline{r},1]^{U}\times \mathbb{R}^{U}.\end{cases}
\end{align}}%
Assumption~\ref{assumption:Lipschitz}~\ref{enum:gradient_Lipschitz} holds with, e.g.,
$\varpi_{1}\coloneqq 4((2+\sqrt{U})\norm{\bm{H}}_{\rm op}^{2} + \norm{\bm{H}^{\T}\bm{y}}) + 2\sqrt{U}\lambda_{r}\underline{r}^{-4} + 2^{-1}\sqrt{U}\lambda_{\theta}M$
and
$\varpi_{2}\coloneqq 4^{-1}M^{2}$,
and thus
a stationary point of~\eqref{eq:MIMO_proposed} can be obtained by Alg.~\ref{alg:proposed} (see Thm.~\ref{theorem:convergence_extension}) 

\subsection{Numerical Experiments}
We evaluated numerical performance of the proposed model~\eqref{eq:MIMO_proposed} and Alg.~\ref{alg:proposed} under the setting\footnote{
  For the problem~\eqref{eq:MIMO_model}, we randomly chose
  (i)
  $\mathsf{s}^{\star} \in \mathsf{D}$;
  (ii)
  $\mathsf{H} \coloneqq \sqrt{\mathsf{R}}\mathsf{G}$
  with
  $\mathsf{G} \in \mathbb{C}^{B\times U}$
  whose entries were sampled from the complex Gaussian distribution
  $\mathbb{C}\mathcal{N}(0,1/U)$,
  and
  a symmetric Toeplitz matrix
  $\mathsf{R} \in \mathbb{R}^{B\times B}$
  whose entries were given by
  $[\mathsf{R}]_{j,k} = 0.5^{\abs{j-k}}$;
  and
  (iii)
  $\mathsf{e} \in \mathbb{C}^{B}$
  whose entries were sampled from
  $\mathbb{C}\mathcal{N}(0,\sigma^{2})$,
  where
  $\sigma^{2}$
  was chosen so that
  $10\log_{10}\frac{1}{\sigma^{2}}\,\mathrm{(dB)}$
  achieved a given signal-to-noise ratio (SNR).
}
in~\cite{Chen19}.
For Alg.~\ref{alg:proposed}, we employed
$(\mu_{n})_{n=1}^{\infty}\coloneqq ((2\eta)^{-1}n^{-1/3})_{n=1}^{\infty}$,
$c = 2^{-13}$,
$\eta = 1$,
and
$(\gamma_{n})_{n=1}^{\infty}$
given by the backtracking algorithm (see Remark~\ref{remark:example}~\ref{enum:ex:stepsize}) with
$\gamma_{\rm initial} \coloneqq 1$
and
$\rho \coloneqq 2^{-1}$.
All experiments were performed by MATLAB on MacBookPro (Apple M3, 16 GB).

We firstly compared convergence performance of Alg.~\ref{alg:proposed} with that of a fairly standard nonsmooth optimization algorithm, called proximal subgradient method
(denoted by ``Sub'')~\cite{Zhu-Zhao-Zhang23}, in a scenario of the problem~\eqref{eq:MIMO_proposed}.
``Sub'' in~\cite{Zhu-Zhao-Zhang23} updates
$\bm{x}_{n+1}\coloneqq \prox{\gamma_{n}\phi}(\bm{x}_{n}-\gamma_{n}\bm{v}_{n}) \in \mathcal{X}$
at
$n$th iteration
with
$\bm{x}_{n} \in \mathcal{X}$,
$\bm{v}_{n} \in \partial (h+g\circ \mathfrak{S})(\bm{x}_{n})$
and
$\gamma_{n} \in \mathbb{R}_{++}$.
For ``Sub'',
we employed two stepsizes:
(i)
$\gamma_{n}\coloneqq (2\varpi_{1}n)^{-1}$,
with
$\varpi_{1}$
given in Sect.~\ref{sec:numerical_formulation},
which guarantees a convergence to a stationary point~\cite[Lem. 4.3 and Thm. 4.1]{Zhu-Zhao-Zhang23};
and
(ii)
heuristic
$\gamma_{n}\coloneqq (2n)^{-1}$\footnote{
  Although the difference between two stepsizes is only a constant factor,
  convergence is not guaranteed for ``Sub'' with the latter heuristic stepsize.
}.
We employed parameters in~\eqref{eq:MIMO_proposed} as
$\lambda_{r}=\lambda_{\theta} = \underline{r} = 0.1$
under
$(U,B,M,\mathrm{SNR})=(128,128,8,20\;\mathrm{(dB)})$.
All algorithms were terminated when running CPU time exceeded
$0.1\;\mathrm{(s)}$.

Fig.~\ref{fig:convergence} shows the averaged values of the cost function over
$100$
trials versus CPU time (s) for Alg.~\ref{alg:proposed} and ``Sub'', where
markers were put at every $100$ iterations.
From Fig.~\ref{fig:convergence}, we observe that Alg.~\ref{alg:proposed} converges much faster than ``Sub'' with stepsize
$\gamma_{n}=(2\varpi_{1}n)^{-1}$
with guaranteed convergence.
Moreover, Alg.~\ref{alg:proposed} also achieves faster convergence speed than ``Sub'' with heuristic stepsize
$\gamma_{n}=(2n)^{-1}$.
This result implies an effective convergence performance of Alg.~\ref{alg:proposed}.

Next, we compared estimation performance of
(i) the proposed model~\eqref{eq:MIMO_proposed}\footnote{
  $\lambda_{r}, \lambda_{\theta}$
  were chosen from
  $\{10^{k}\mid k=-6,-5,\ldots,0\}$
  and
  $\underline{r} = 0.1$.
}, solved by Alg.~\ref{alg:proposed},
with that of three models: (ii) ``LMMSE'', e.g.,~\cite{Yang-Hanzo15}, (iii)
``Modulus''\footnote{
  We applied a projected gradient method~\cite{Chen19} to
  ``Modulus''.
}~\cite{Chen19},
and (iv) ``SOAV''\footnote{
  $\lambda$
  in~\eqref{eq:MIMO_problem} was chosen from
  $\{10^{k}\mid k=-6,-5,\ldots,1\}$.
  We applied a primal-dual splitting algorithm~\cite{Condat13} to ``SOAV''.
}~\cite{Hayakawa-Hayashi20}
(see just after~\eqref{eq:MIMO_problem} for (ii)-(iv) in detail).
All algorithms except ``LMMSE'' were terminated when
running CPU time exceeded
$3\;\mathrm{(s)}$
or
$\norm{\bm{x}_{n+1}-\bm{x}_{n}} \leq 10^{-5}$
held with a sequence
$(\bm{x}_{n})_{n=1}^{\infty}$
generated by algorithms.
As a performance criterion, we employed averaged bit error rate (BER)\footnote{
  BER was computed by a MATLAB code with
  'pskdemod' and 'biterr'.
} of the final estimate over
$100$
trials.
The lower averaged BER indicates better estimate performance.

Figs.~\ref{fig:BER1} and~\ref{fig:BER2} demonstrate the averaged BER of each algorithm versus SNR under respectively two settings
$B=U$
and
$B=\frac{3}{4}U$,
where the averaged BER
$0$
was replaced by the machine epsilon.
We note that the estimation task~\eqref{eq:MIMO_model} becomes challenging as the ratio
$B/U$
becomes small.
From Figs.~\ref{fig:BER1} and~\ref{fig:BER2}, we observe that (i) every model achieves better performance than ``LMMSE'';
(ii) the proposed model in~\eqref{eq:MIMO_proposed} outperforms the others for all SNRs and ratios
$B/U$.
In particular, from Fig.~\ref{fig:BER2},  we see that the proposed model in~\eqref{eq:MIMO_proposed} overwhelms the others in the challenging setting with
$B=\frac{3}{4}U$.
As we expected, these results show a great potential of the proposed model in~\eqref{eq:MIMO_proposed} and Alg.~\ref{alg:proposed} for MU-MIMO signal detection.

\section{Conclusion}
We proposed a proximal variable smoothing for nonsmooth minimization involving weakly convex composite in Problem~\ref{problem:origin}.
The proposed algorithm consists of two steps: (i) a time-varying forward step with the gradient of a smoothed surrogate function designed with the Moreau envelope; (ii) the backward step with the proximity operator.
We also presented an asymptotic convergence analysis of the proposed algorithm.
Numerical experiments in a scenario of MU-MIMO demonstrate the effectiveness of (i) the proposed algorithm and (ii) a new formulation in~\eqref{eq:MIMO_proposed} of MU-MIMO by Problem~\ref{problem:origin}.

\bibliographystyle{IEEEtran}
\bibliography{refs}

\begin{thebibliography}{10}
\providecommand{\url}[1]{#1}
\csname url@samestyle\endcsname
\providecommand{\newblock}{\relax}
\providecommand{\bibinfo}[2]{#2}
\providecommand{\BIBentrySTDinterwordspacing}{\spaceskip=0pt\relax}
\providecommand{\BIBentryALTinterwordstretchfactor}{4}
\providecommand{\BIBentryALTinterwordspacing}{\spaceskip=\fontdimen2\font plus
\BIBentryALTinterwordstretchfactor\fontdimen3\font minus \fontdimen4\font\relax}
\providecommand{\BIBforeignlanguage}[2]{{%
\expandafter\ifx\csname l@#1\endcsname\relax
\typeout{** WARNING: IEEEtran.bst: No hyphenation pattern has been}%
\typeout{** loaded for the language `#1'. Using the pattern for}%
\typeout{** the default language instead.}%
\else
\language=\csname l@#1\endcsname
\fi
#2}}
\providecommand{\BIBdecl}{\relax}
\BIBdecl

\bibitem{Bauschke-Combettes17}
H.~H. Bauschke and P.~L. Combettes, \emph{Convex analysis and monotone operator theory in Hilbert spaces}, 2nd~ed.\hskip 1em plus 0.5em minus 0.4em\relax Springer, 2017.

\bibitem{Zhang10}
C.-H. Zhang, ``{Nearly unbiased variable selection under minimax concave penalty},'' \emph{The Annals of Statistics}, vol.~38, no.~2, pp. 894 -- 942, 2010.

\bibitem{Fan-Li01}
J.~Fan and R.~Li, ``\BIBforeignlanguage{English (US)}{Variable selection via nonconcave penalized likelihood and its oracle properties},'' \emph{\BIBforeignlanguage{English (US)}{Journal of the American Statistical Association}}, vol.~96, no. 456, pp. 1348--1360, 2001.

\bibitem{Li-So-Ma20}
J.~Li, A.~M.~C. So, and W.~K. Ma, ``Understanding notions of stationarity in nonsmooth optimization: A guided tour of various constructions of subdifferential for nonsmooth functions,'' \emph{IEEE Signal Processing Magazine}, vol.~37, no.~5, pp. 18--31, 2020.

\bibitem{Gandy-Recht-Yamada11}
S.~Gandy, B.~Recht, and I.~Yamada, ``Tensor completion and low-n-rank tensor recovery via convex optimization,'' \emph{Inverse Problems}, vol.~27, no.~2, p. 025010, 2011.

\bibitem{Yin-Parekh-Selesnick19}
L.~Yin, A.~Parekh, and I.~Selesnick, ``Stable principal component pursuit via convex analysis,'' \emph{IEEE Transactions on Signal Processing}, vol.~67, no.~10, pp. 2595--2607, 2019.

\bibitem{Kuroda-Kitahara22}
H.~Kuroda and D.~Kitahara, ``Block-sparse recovery with optimal block partition,'' \emph{IEEE Transactions on Signal Processing}, vol.~70, pp. 1506--1520, 2022.

\bibitem{Combettes-Pesquet21}
P.~L. Combettes and J.-C. Pesquet, ``Fixed point strategies in data science,'' \emph{IEEE Transactions on Signal Processing}, vol.~69, pp. 3878--3905, 2021.

\bibitem{Zhen-Ma-Xue24}
Z.~Zheng, S.~Ma, and L.~Xue, ``A new inexact proximal linear algorithm with adaptive stopping criteria for robust phase retrieval,'' \emph{IEEE Transactions on Signal Processing}, vol.~72, pp. 1081--1093, 2024.

\bibitem{Duchi-Ruan18}
J.~C. Duchi and F.~Ruan, ``{Solving (most) of a set of quadratic equalities: composite optimization for robust phase retrieval},'' \emph{Information and Inference: A Journal of the IMA}, vol.~8, no.~3, pp. 471--529, 2018.

\bibitem{Charisopoulos-Chen-Davis-Diaz-Ding-Drusvyatskiy21}
V.~Charisopoulos, Y.~Chen, D.~Davis, M.~D{\'i}az, L.~Ding, and D.~Drusvyatskiy, ``Low-rank matrix recovery with composite optimization: Good conditioning and rapid convergence,'' \emph{Foundations of Computational Mathematics}, vol.~21, no.~6, pp. 1505--1593, 2021.

\bibitem{Wang-So-Zoubir23}
Z.-Y. Wang, H.~C. So, and A.~M. Zoubir, ``Robust low-rank matrix recovery via hybrid ordinary-{W}elsch function,'' \emph{IEEE Transactions on Signal Processing}, vol.~71, pp. 2548--2563, 2023.

\bibitem{Charisopoulos-Davis-Dias-Drusvyatskiy20}
V.~Charisopoulos, D.~Davis, M.~D^^c3^^adaz, and D.~Drusvyatskiy, ``{Composite optimization for robust rank one bilinear sensing},'' \emph{Information and Inference: A Journal of the IMA}, vol.~10, no.~2, pp. 333--396, 2020.

\bibitem{Lewis-Wright16}
A.~S. Lewis and S.~J. Wright, ``A proximal method for composite minimization,'' \emph{Mathematical Programming}, vol. 158, no.~1, pp. 501--546, 2016.

\bibitem{Drusvyatskiy-Paquette19}
D.~Drusvyatskiy and C.~Paquette, ``Efficiency of minimizing compositions of convex functions and smooth maps,'' \emph{Mathematical Programming}, vol. 178, no.~1, pp. 503--558, 2019.

\bibitem{Davis-Yin17}
D.~Davis and W.~Yin, ``{A Three-Operator Splitting Scheme and its Optimization Applications},'' \emph{Set-Valued and Variational Analysis}, vol.~25, no.~4, pp. 829--858, 2017.

\bibitem{Zhao-Cevher18}
R.~Zhao and V.~Cevher, ``Stochastic three-composite convex minimization with a linear operator,'' in \emph{AISTATS}, vol.~84, 2018, pp. 765--774.

\bibitem{Condat-Kitahara-Contreras-Hirabayashi23}
L.~Condat, D.~Kitahara, A.~Contreras, and A.~Hirabayashi, ``Proximal splitting algorithms for convex optimization: a tour of recent advances, with new twists,'' \emph{SIAM Review}, vol.~65, no.~2, pp. 375--435, 2023.

\bibitem{Bohm-Weight21}
A.~B{\"o}hm and S.~J. Wright, ``Variable smoothing for weakly convex composite functions,'' \emph{Journal of Optimization Theory and Applications}, vol. 188, no.~3, pp. 628--649, 2021.

\bibitem{Kume-Yamada24}
K.~Kume and I.~Yamada, ``A variable smoothing for nonconvexly constrained nonsmooth optimization with application to sparse spectral clustering,'' in \emph{IEEE ICASSP}, 2024, pp. 9296--9300.

\bibitem{Liu-Xia24}
Y.~Liu and F.~Xia, ``{Proximal variable smoothing method for three-composite nonconvex nonsmooth minimization with a linear operator},'' \emph{Numerical Algorithms}, vol.~96, no.~1, pp. 237--266, 2024.

\bibitem{Kume-Yamada24B}
K.~Kume and I.~Yamada, ``A variable smoothing for weakly convex composite minimization with nonconvex constraint,'' \emph{arXiv (2412.04225)}, pp. 1--39, 2024.

\bibitem{Yang-Hanzo15}
S.~Yang and L.~Hanzo, ``Fifty years of {MIMO} detection: The road to large-scale {MIMO}s,'' \emph{IEEE Communications Surveys \& Tutorials}, vol.~17, no.~4, pp. 1941--1988, 2015.

\bibitem{Chen19}
J.~C. Chen, ``Computationally efficient data detection algorithm for massive {MU-MIMO} systems using {PSK} modulations,'' \emph{IEEE Communications Letters}, vol.~23, no.~6, pp. 983--986, 2019.

\bibitem{Hayakawa19}
R.~Hayakawa and K.~Hayashi, ``Discrete-valued vector reconstruction by optimization with sum of sparse regularizers,'' in \emph{EUSIPCO}.\hskip 1em plus 0.5em minus 0.4em\relax EURASIP, 2019, pp. 1--5.

\bibitem{Hayakawa-Hayashi20}
------, ``Asymptotic performance of discrete-valued vector reconstruction via box-constrained optimization with sum of l1 regularizers,'' \emph{IEEE Transactions on Signal Processing}, vol.~68, pp. 4320--4335, 2020.

\bibitem{Rockafellar-Wets98}
R.~Rockafellar and R.~J.-B. Wets, \emph{Variational Analysis}, 3rd~ed.\hskip 1em plus 0.5em minus 0.4em\relax Springer Verlag, 2010.

\bibitem{Yamada-Yukawa-Yamagishi11}
I.~Yamada, M.~Yukawa, and M.~Yamagishi, ``Minimizing the {M}oreau envelope of nonsmooth convex functions over the fixed point set of certain quasi-nonexpansive mappings,'' in \emph{Fixed-Point Algorithms for Inverse Problems in Science and Engineering}, H.~H. Bauschke, R.~S. Burachik, P.~L. Combettes, V.~Elser, D.~R. Luke, and H.~Wolkowicz, Eds.\hskip 1em plus 0.5em minus 0.4em\relax Springer New York, 2011, pp. 345--390.

\bibitem{Abe-Yamagishi-Yamada20}
J.~Abe, M.~Yamagishi, and I.~Yamada, ``Linearly involved generalized {M}oreau enhanced models and their proximal splitting algorithm under overall convexity condition,'' \emph{Inverse Problems}, vol.~36, no.~3, 2020.

\bibitem{Bauschke-Moursi-Wang21}
H.~H. Bauschke, W.~M. Moursi, and X.~Wang, ``Generalized monotone operators and their averaged resolvents,'' \emph{Mathematical Programming}, vol. 189, no.~1, pp. 55--74, 2021.

\bibitem{prox_repository}
\BIBentryALTinterwordspacing
G.~Chierchia, E.~Chouzenoux, P.~L. Combettes, and J.-C. Pesquet, ``{The proximity operator repository}.'' [Online]. Available: \url{http://proximity-operator.net/}
\BIBentrySTDinterwordspacing

\bibitem{Hoheisel-Laborde-Oberman20}
T.~Hoheisel, M.~Laborde, and A.~Oberman, ``A regularization interpretation of the proximal point method for weakly convex functions,'' \emph{Journal of Dynamics and Games}, vol.~7, no.~1, pp. 79--96, 2020.

\bibitem{Levin-Kileel-Boumal23}
E.~Levin, J.~Kileel, and N.~Boumal, ``Finding stationary points on bounded-rank matrices: a geometric hurdle and a smooth remedy,'' \emph{Mathematical Programming}, vol. 199, no.~1, pp. 831--864, 2023.

\bibitem{Beck17}
A.~Beck, \emph{First-Order Methods in Optimization}.\hskip 1em plus 0.5em minus 0.4em\relax SIAM, 2017.

\bibitem{Nagahara15}
M.~Nagahara, ``Discrete signal reconstruction by sum of absolute values,'' \emph{IEEE Signal Processing Letters}, vol.~22, no.~10, pp. 1575--1579, 2015.

\bibitem{Zhu-Zhao-Zhang23}
D.~Zhu, L.~Zhao, and S.~Zhang, ``A unified analysis for the subgradient methods minimizing composite nonconvex, nonsmooth and non-{L}ipschitz functions,'' \emph{arXiv (2308.16362)}, pp. 1--32, 2023.

\bibitem{Condat13}
L.~Condat, ``A primal^^e2^^80^^93dual splitting method for convex optimization involving {L}ipschitzian, proximable and linear composite terms,'' \emph{Journal of Optimization Theory and Applications}, vol. 158, no.~2, pp. 460--479, 2013.

\end{thebibliography}

\end{document}